\documentclass[reqno,11pt]{article}

\usepackage{amsmath,latexsym,amssymb}

\def\sqr#1#2{{\vcenter{\hrule height.#2pt
        \hbox{\vrule width.#2pt height#1pt \kern#1pt
                \vrule width.#2pt}
        \hrule height.#2pt}}}

\newtheorem{Theorem}{\sc Theorem}[section]
\newtheorem{Lemma}[Theorem]{\sc Lemma}
\newtheorem{Corollary}[Theorem]{\sc Corollary}
\newtheorem{Proposition}[Theorem]{\sc Proposition}

\newtheorem{Remark}[Theorem]{\sc Remark}

\newtheorem{Definition}[Theorem]{\sc Definition}

\makeatletter
\def\im{\mathop{\operator@font im}\nolimits}
\def\red{{\mathop{\operator@font red}\nolimits}}
\def\arg{\mathop{\operator@font arg}\nolimits}

\def\frk{\mathop{\mbox{f-rk}}\nolimits}
\def\pole{\mathop{\operator@font pole}\nolimits}
\def\tor{{\mathop{\operator@font tor}\nolimits}}
\def\Cst{\mathop{\operator@font Cst}\nolimits}
\def\Newton{\mathop{\operator@font Newton}\nolimits}
\def\Hom{\mathop{\operator@font Hom}\nolimits}
\def\Ext{\mathop{\operator@font Ext}\nolimits}
\def\End{\mathop{\operator@font End}\nolimits}
\def\Ass{\mathop{\operator@font Ass}\nolimits}
\def\rk{\mathop{\operator@font rk}\nolimits}
\def\ker{\mathop{\operator@font ker}\nolimits}
\def\coker{\mathop{\operator@font coker}\nolimits}

\def\m{\mathop{\operator@font m}\nolimits}

\newcommand{\htt}{\mathop{\operator@font ht}\nolimits}
\def\Ann{\mathop{\operator@font Ann}\nolimits}
\def\MCM{\mathop{\operator@font MCM}\nolimits}
\def\FID{\mathop{\operator@font FID}\nolimits}
\def\MC{\mathop{\operator@font MC}\nolimits}
\def\CEC{\mathop{\operator@font CEC}\nolimits}
\def\gr{\mathop{\operator@font gr}\nolimits}
\def\Tor{\mathop{\operator@font Tor}\nolimits}
\def\ext{\mathop{\operator@font ext}\nolimits}

\makeatother

\begin{document}

\baselineskip=13pt

\pagestyle{empty}

\ \vspace{1.7in}

\noindent {\LARGE\bf Gorenstein rings call the tune}

\vspace{.25in}

\noindent \ Anne-Marie Simon, \  Service d'Alg\`{e}bre, C.P. 211, Universit\'{e}
Libre de Bruxelles, 1050 Bruxelles, Belgique. {\it E-mail}: {\ amsimon@ulb.ac.be}

\bigskip

\noindent  \ Jan R. Strooker, \  Mathematisch Instituut,
Universiteit Utrecht, Postbus 80010, 3508 TA Utrecht, Nederland. {\it E-mail}: {\
strooker@math.uu.nl

\vspace{2.4cm}

\section{Introduction \hfill\break}

The `Homological Conjectures' in local algebra date back to Serre's 
beautiful 1957-58 course at the Coll\`{e}ge de France \cite
{Se}. How to count the 
multiplicities of components in an intersection of two algebraic varieties, 
was reduced to questions in local algebra involving the celebrated Tor-formula.
Several important results were proven, other questions remained. These led
to related but also different conjectures by other mathematicians, notably
M. Auslander, Bass \cite{Ba} Vasconcelos \cite{Va}, Peskine-Szpiro \cite{PS}, and Hochster \cite{Ho 1}, \cite{Ho 3}. 
Some of these were proved for all noetherian local rings, others up to a certain dimension or in equal characteristic; others again, remain virtually untouched.  A masterly survey of results and of the many interconnections between these 
Homological Conjectures was presented in \cite{Ho 2}. This was updated in \cite{Ro 1}. 
Progress depending on De Jong's theory of alterations is discussed in \cite{Be}, \cite{Ro 2}. Finally, a substantial part of these conjectures is treated more leisurely in the monographs \cite{St} and \cite{BH}. 
	
Our aim in this note is far less ambitious. We focus on two of 
Hochster's conjectures: the Monomial Conjecture (MC) and the Canonical Element
Conjecture (CEC). We point out that these are equivalent with certain statements about Gorenstein rings, in a way to be made precise as we go along. Thus 
Gorenstein rings `call the tune': it is their module theory, or if one likes, representation theory, which controls MC and CEC for all noetherian local rings. To us, this suggests that the Gorenstein property, besides the several `symmetries' and `dualities' known, harbours further ones awaiting discovery.

It is not customary to write a survey of material already in print. However,
we think that a brief presentation of selected results from a series of papers 
\cite{St2}, \cite{SS 1}, \cite{SS 2}, \cite{SS 3}, and \cite{SS 4}, highlighting the Gorenstein connexion, will help the reader to understand our point of view and perhaps continue this line of research. No proofs are given, rather remarks about ideas and results used along the way. We are grateful to the editors for this opportunity.  

\bigskip

\section{The conjectures \hfill\break}

We use the word `ring' for unital, commutative, noetherian local rings,
since no others shall be needed, and all modules will be finitely generated
unless otherwise stated. As a standard notation we use $(A, \mathfrak m, k)$ for the ring $A$, its maximal ideal and its residue class field. Let $A$ be a $d$-dimensional ring, and $x_{1}, ..., x_{d} = \bf x$ a system of parameters (sop). Let $t$ be a positive integer. Is the product $x_{1}^{t}...x_{d}^{t}$ in the ideal generated by $x_{1}^{t+1}, ..., x_{d}^{t+1}$?

The ring is said to satisfy MC provided for every choice of $\bf x$ and of $t$,
one has $x_{1}^{t}...x_{d}^{t} \not\in (x_{1}^{t+1}, ..., x_{d}^{t+1})$  \cite[Conj. 2]{Ho 1}.
It is often advantageous to translate this to saying that the maps
$x_{1}^{t}...x_{d}^{t}:	A/(x_{1}, ..., x_{d}) \rightarrow A/(x_{1}^{t+1}, ...,
x_{d}^{t+1})$ are never null. Since the direct limit of this directed system 
is the $d$-th local cohomology module $H_{\mathfrak m}^{d}(A)$, we can also state
\begin{Definition}{\rm A ring $A$ is said to satisfy the Monomial Conjecture if
the induced map $\mu_{\bf x}(A): A/(x_{1}, ..., x_{d}) \rightarrow H_{\mathfrak m}^{d}(A)$ is nonnull for every sop $\bf x$.}
\end{Definition}

To prepare for the Canonical Element Conjecture, let $M$ be a finitely 
generated $A$-module and consider the natural maps $\theta^{i}_{A}(M):
\Ext^{i}_{A}(k,M) \rightarrow H_{\mathfrak m}^{i}(k,M)$ for $i \geq 0$. 
This Ext is a finite dimensional vector space over $k$, whose dimension is often called the $i$-th Bass number $\mu^{i}_{A}(M)$ of $M$; we call the $k$-dimension of the image of $\theta^{i}_{A}(M)$ the reduced Bass number $\nu^{i}_{A}(M)$ of $M$. Thus $\nu^{i}_{A}(M) \leq \mu^{i}_{A}(M)$ and both are nonnegative integers.

Let $\bf F$ be a resolution of $k$ in terms of free $A$-modules. Let $\bf x$ again be an sop, and ${\bf K}({\bf x},A)$ its Koszul complex. Then the surjection
$A/(x_{1}, ..., x_{d}) \rightarrow k$ lifts to a map of complexes 
$\phi: {\bf K}({\bf x},A) \rightarrow \bf {F}$. 
\begin{Definition}{\rm A ring $A$ is said to satisfy the Canonical Element
Conjecture if the map $\phi_{d}$ is nonnull for every lifting belonging to every sop $\bf x$.}    
\end{Definition} 
 
Notice that the source of $\phi_{d}$ is just a copy of $A$. This was Hochster's original approach \cite[Def. 2.1]{Ho 3}. He and others have put forward several equivalent formulations of CEC. One of these led to \cite[Ths. 3.7 and 4.3]{Ho 3}, see also \cite[5.1]{SS 3}:
\begin{Theorem} Write $S \subset F_{d}$ for the $d$-th syzygy of $k$.
Then the following are equivalent:
\begin{enumerate}
\item[(i)] The ring $A$ satisfies the Canonical Element Conjecture; 
\item[(ii)]$\nu^{d}_{A}(S) > 0$; 
\item[(iii)]$\nu^{d}_{A}(M) > 0$ for some $M$, not necessarily finitely generated.
\end{enumerate}
In case $A$ possesses a canonical module $K$, this is also equivalent with 
\begin{enumerate}
\item[(iv)]$\nu^{d}_{A}(K) > 0$.
\end{enumerate} 
\end{Theorem}
	Here $K$ is called a canonical module of $A$ when ${K}^{\vee}$ = $H_{\mathfrak m}^{d}(A)$, where ${-}^{\vee}$ stands for the Matlis dual, e.g.
\cite[1.6]{SS 3}. Both MC and CEC need only be proved for rings which are complete in 
their $\mathfrak m$-adic topology, and a ring which satisfies CEC satsifies MC. In case the sop $\bf x$ forms a regular sequence, i.e. $A$ is a Cohen-Macaulay ring, both conjectures are true. Better even, if there is an $A$-module $M$, not necessarily finitely generated, such that $x_{1}, ..., x_{d}$ form a regular sequence on $M$ and $M \neq (x_{1}, ..., x_{d})M$, then both remain true for $\hat{M}$, the $\mathfrak m$-adic completion of $M$, and there  hold for all sop's $\bf x$ \cite[Th. 9.1.1]{St}. Such, possibly `Big', Cohen-Macaulay modules, provide one of the proofs of the conjectures in equal characteristic. All this, and more, is in the papers \cite{Ho 1} and \cite{Ho 3}. Thanks to Heitmann's recent advance in dimension 3, \cite{He} and \cite{Ro 3}, MC and CEC are known in mixed characteristic for $d \leq 3$ but not beyond. It should be mentioned that the two conjectures are closely interwoven with others, in particular Hochster's Direct Summand Conjecture, but we don't deal with these.	

Here we wish to emphasize that it is the connexion with local cohomology and
then the description in terms of reduced Bass numbers, which shall allow us 
to tie up the conjectures with properties of Gorenstein rings.
\section{Why Gorenstein?}
Since we need only prove MC for rings complete in their maximal ideal
topology, let $A$ be such a ring. According to 
Cohen, one can write $A = S/\mathfrak b$ where b is an ideal in the regular ring
$S$. By judiciously dividing out a maximal regular sequence $z_{1}, ..., z_{n}$
in $\mathfrak b$, one can prove \cite[Prop. 1]{St2}, \cite[5.5]{SS 3}
\begin{Lemma} There exists a $d$-dimensional complete intersection ring $R$ 
with $d = \dim A$ which maps onto $A = R/\mathfrak a$. Moreover, sop's in $R$ map to sop's in $A$, and any sop $\bf x$ in $A$ can be lifted to an sop $\bf y$ in $R$.
\end{Lemma}

In the situation of this Lemma, $\mu_{\bf x}(A) = \mu_{\bf y}(R) \otimes_{R}A$   and $\mu_{\bf x}(A) = 0$ if and only if its Matlis dual $\mu_{\bf x}(A)^{\vee}$
is.
Working out this dual, using the adjointness of tensor product and Hom
and a few facts from local Grothendieck duality regarding the Gorenstein ring $R$, one obtains \cite[Prop 2]{St2}:
\begin{Proposition} In the situation just described, $\mu_{\bf x}(A) = 0$
if and only if $\Ann_{R}\mathfrak a \subset (y_{1}, ..., y_{d})$. 
\end{Proposition}
	
\begin{Definition}{\rm We say that a Gorenstein ring $R$ and a ring $A$ are in our favourite position if $A = R/\mathfrak a$ where $\mathfrak a$ in a nonnull ideal in $R$ consisting of zero divisors. We standardly put $\mathfrak b = \Ann_{R} \mathfrak a$.}
\end{Definition} 

With the help of the above and further reasoning, one proves 
\begin{Theorem} Let $R$ and $A$ be in our favourite position. Then $A$ satisfies
$\MC$ if and only if the ideal $\mathfrak b$ is not contained in any parameter ideal of $R$.
\end{Theorem}

It is well known that in a Gorenstein ring the ideals $\mathfrak b = \Ann \mathfrak a$ where $\mathfrak a$ is a nonnull ideal consisting of zero divisors,
are just the unmixed nonnull ideals consisting of zero divisors. So one has
\begin{Corollary} Every ring satisfies the Monomial Conjecture if and only 
if in no Gorenstein ring is an unmixed nonnull ideal of zero divisors contained in a parameter ideal.
\end{Corollary}
\begin{Remark}{\rm In our favourite situation, if $R$ is equicharacteristic, so is the ring $A$ which therefore satisfies MC. It follows that the above statement is true for every equicharacteristic Gorenstein ring. It only remains 
an open question in mixed characteristic for (Heitmann) dimension $> 3$}.
\end{Remark}
\begin{Remark}{\rm Furthermore, in the corollary we may replace the word Gorenstein by complete intersection. At face value, the statement becomes stronger in one direction, weaker in the other. We have not managed so far to utilize properties of complete intersections which are not shared by all Gorenstein rings to prove MC or CEC in cases yet unknown. Also our feeling is that the dualities and symmetries pertaining to Gorenstein rings and their module theory are just `right'. So here, and in following sections, we shall stick to Gorenstein and leave it to the reader to realize that if one can obtain
certain results for complete intersections, this is enough}.
\end{Remark} 

For details and further observations we refer to \cite{St2}.
\section{Auslander-Buchweitz theory}
The theory which M. Auslander and R.-O. Buchweitz developed in \cite{AB}, is
quite abstract and covers various special cases, of which we treat 
only one: (finitely generated) modules over (local) Gorenstein rings. For the 
sake of brevity, we recall this in a manner which is not quite that of \cite{AB}, but was worked out later by various authors and is systematized in the first four sections of \cite{SS 3}. In this section we fix a Gorenstein ring $R$ and only consider modules and several invariants over this. We therefore drop the $R$ from notation. 

Let $M$ be a module over the Gorenstein ring $R$ and $C$ a maximal Cohen-Macaulay module (MCM). An $R$-homomorphism $f: C \rightarrow M$ is called an MCM-precover of $M$ if, for every MCM-module $D$, the induced map
$\Hom(D,C) \rightarrow \Hom(D,M)$ is surjective. In other words, every map
from $D$ to $M$ can be factored through $C$. If, in addition, the map $f$ is right minimal, then it is called a cover. Here right minimal means that if $f \circ g = f$ for any endomorphism $g$ of $C$, then $g$ is an automorphism. 

	Dually, one defines a homomorphism $h: M \rightarrow J$ with $J$ a
module of finite injective dimension (FID) to be an FID-preenvelope resp. 
-envelope of $M$. Very briefly, the basic result of \cite{AB}
can in our case be summarized \cite[section 3]{SS 3} as
\begin{Theorem} Over a Gorenstein ring, there exist $\MCM$-covers and $\FID$-envelopes. In these cases $f$ is a surjection and $h$ an injection with an $\FID$-kernel resp. $\MCM$-cokernel. Covers and envelopes are determined up to isomorphisms over $M$.
\end{Theorem} 
	
Though we shall stick with our Gorenstein environment, it should be pointed out that `Auslander-Buchweitz contexts' turn out to be relevant in different
areas of algebra \cite{Ha}, beyond what was originally envisaged by these authors.

By the $\frk$ of a module $X$, we mean the maximal rank of a free direct
summand of $X$. For a module $M$, by the above theorem, the $\frk$ of its MCM-cover $C$ and its FID-envelope $J$ are well defined. In \cite[p. 8]{AB} it is suggested to examine these. The first one is called Auslander's delta invariant $\delta(M)$ and has been investigated by several of his last students and others. The one on the FID-side so far has hardly been looked at. We hope to convince you 
that this is unjustified, and that these invariants are most profitably studied in tandem.   
\begin{Proposition} Let $J$ be an $\FID$-envelope of $M$. Then $\frk(J) =
\nu^{d}(M)$.
\end{Proposition} 

	This result, \cite[Prop. 6]{SS 2}, \cite[Prop. 3.10]{SS 3}, connects the second invariant with local cohomology. What about the $\delta$-invariant? Well,
\cite[Prop. 7]{SS 2} and \cite[Prop. 4.1 (iv)]{SS 3}:
\begin{Proposition} Let $p: R^{t} \rightarrow M$ be a surjection of a free module onto $M$. This induces a map $\Ext^{d}(k,R^{t}) \rightarrow 
\Ext^{d}(k,M)$. The first $\Ext$ being simply $k^{t}$, its image is a 
finite dimensional vector space $V$ over $k$. Then $\delta(M) = \dim V$.
\end{Proposition}

	While this description of $\delta$ does not involve local 
cohomology, it does bring into play the $d$-th $\Ext$ whose relation 
with $d$-th local cohomology is crucial for the reduced Bass number. 
We shall exploit this affinity in the next section, returning to our main theme. 
There are many interactions between the two invariants, which suggest 
a formal kind of duality. Perhaps more accurately one should speak of 
orthogonality, stemming from the fact that if $X$ is an MCM and $Y$ an 
FID, then $\Ext^{i}(X,Y) = 0$ for $i > 0$. Of course, the yoga becomes
even more powerful when one realizes that over a Gorenstein ring,
FID-modules are exactly the ones of finite projective dimension \cite[Th. 10.1.9]{St}.
For our purposes, we collect a few useful facts from this wondrous world:
\begin{Lemma} For any module $M$, there exists a map $s: M \rightarrow R$
such that $\Ext^{d}(k,s)$ is nonnull if and only if $\mathfrak \nu^{d}(M)
> 0$.
\end{Lemma}
Indeed, bearing in mind that the canonical module of the Gorenstein ring $R$
is just $R$ itself, one concludes this from \cite[Cor. 3.11]{SS 3}. A more precise statement, which can be regarded as a counterpart on the envelope side
of Proposition 4.3, can be found in \cite[Th. 3]{SS 2} or \cite[Cor. 3.12]{SS 3}. 

The next lemma is \cite[Prop. 4.1 (iii), Prop. 4.5]{SS 3}.
\begin{Lemma} For any module $M$, $\delta(M) \leq \beta(M)$, where 
$\beta(M)$ means the minimal number of generators of the module $M$. 
For $\FID$-modules this is actually an equality. Also $\delta(M) = 0$
precisely when $M$ is a homomorphic image of an $\MCM$-module without a free direct summand. 
\end{Lemma}  

\section{Harvesting from favourite position}
Here is one more reason why our favourite position is so effective.
In the set-up of Definition 3.3, the ideal $\mathfrak b$, being annihilated by
$\mathfrak a$, is also an $A$-module. In fact, it is a canonical 
module of the ring $A$ \cite[6.9 and 1.6]{SS 3}.
\begin{Theorem} In our favourite position consider the statements
\begin{enumerate}
\item [(i)]	The ring $A$ satisfies the Canonical Element Conjecture;
\item [(ii)] 	$\nu^{d}_{R}(\mathfrak b) > 0$;
\item [(iii)]	$\delta(R/\mathfrak b) = 0$;
\item [(iv)]	$\mathfrak b$ is not contained in any $\FID$-ideal;
\item [(v)]	$\mathfrak b$ is not contained in any parameter ideal;
\item [(vi)]	The ring $A$ satisfies the Monomial Conjecture;
\end{enumerate}
where items (ii) - (v) all refer to ideals of the Gorenstein ring $R$.
Then (i) $\Rightarrow$ (ii) $\Leftrightarrow$ (iii) $\Rightarrow$ (iv) $\Rightarrow$ (v) $\Leftrightarrow$ (vi).
\end{Theorem}
Sketch of proof: \textit{(i)} $\Rightarrow$ \textit{(ii)}: Theorem 2.3 tells us that \textit{(i)} and
$\nu^{d}_{A}(\mathfrak b) > 0$ are equivalent statements. Since the map $\theta^{d}_{A}(\mathfrak b)$ (section 2) factors through $\theta^{d}_{R}(\mathfrak b)$, one has \textit{(ii)}.
\textit{(ii)} $\Leftrightarrow$ \textit{(iii)}: The embedding $\mathfrak b \rightarrow R$ yields an exact sequence $\Ext^{d}_{R}(k,\mathfrak b) \rightarrow \Ext^{d}_{R}(k,R) \rightarrow \Ext^{d}_{R}(k,R/\mathfrak b)$. Since the middle Ext is just a single copy of $k$ for the Gorenstein ring $R$, and any map $\mathfrak b \rightarrow R$ lands in $\Ann \mathfrak a= \mathfrak b$, Lemma 4.4 and Prop. 4.3 allow one to conclude. \textit{(iii)} $\Rightarrow$ \textit{(iv)}: If $\mathfrak c$ in an $\FID$-ideal containing $\mathfrak b$, then $\delta(R/\mathfrak c) = 1$ by Lemma 4.5. Which contradicts the last statement of this lemma, because $R/\mathfrak c$ is a homomorphic image of $R/\mathfrak b$. \textit{(iv)} $\Rightarrow$ \textit{(v)}: In the Gorenstein ring $R$, every sop is a regular sequence which generates an ideal of finite projective dimension. This then also has $\FID$. \textit{(v)} $\Leftrightarrow$ \textit{(vi)}: Our Theorem 3.4.

	We now recall one of the principal results in \cite{Ho 3}, its Theorem 2.8.
\begin{Theorem} If the Monomial Conjecture is true for all rings, then all
rings satisfy the Canonical Element Conjecture.
\end{Theorem}
Well, this is not quite what is stated there, but it is if one takes into account the equivalence of Direct Summand Conjecture and Monomial Conjecture from the earlier paper \cite[Th. 1]{Ho 1}.

In a more precise statement of Hochster's Theorem 5.2 one may replace  
`all rings' both times by `all rings of a particular dimension and a fixed
residue characteristic'. To a ring ($A,\mathfrak m,k$) we attach the pair of integers ($u,v$) when $u$= char $A$ and $v$ = char $k$. In the cases (0,0) and ($p,p$), $p$ a prime, we speak of equal characteristic, while (0,$p$) is mixed.
For fixed $u$ and $v$ we speak of a particular characteristic.
In our favourite position, if $R$ is equicharacteristic, so is $A$ of the same 
particular characteristic. Now contemplate the statements and implications in
Theorem 5.1, bearing in mind that CEC is known in equal characteristic and 
in the other case for $d \leq 3$. 
\begin{Corollary} All rings of a particular characteristic and dimension satisfy $\CEC$ (or $\MC$, for that matter) if and only if statements (ii) - (v) hold for all nonnull unmixed ideals of zero divisors $\mathfrak b$ in all such Gorenstein rings. 
\end{Corollary}

This gives new information about Gorenstein rings of equal characteristic and 
for mixed characteristic in dimension $\leq 3$. We believe that this should be
exploited and further structure will come to light. In the remaining cases, 
the burden of the Conjectures is borne by the Gorenstein specimens. 

Using several of these results one can prove \cite[Cor. 7.2]{SS 3}.
\begin{Proposition} The Canonical Element Conjecture is equivalent with the
following statement. For any system of parameters in a Gorenstein ring
the kernels of all maps $\phi_{d}$ of Definition 2.2 are contained in the 
nilradical of the ring.
\end{Proposition}
This again suggests some kind of symmetry, at least for the minimal primes in a Gorenstein ring. For further results and refinements see \cite{SS 2} and \cite{SS 3}. 

We finish this section with a tribute to the late Maurice Auslander. Once again ideas of his, this time in collaboration with Buchweitz, which started out as a kind of abstract `Spielerei', revealed themselves over the years as relevant to more concrete problems investigated by others. 
\section{Back to linear algebra: stiffness}  
Let $A$ be a ring and $f: A^{m} \rightarrow A^{n}$ a homomorphism between two free $A$-modules. By choosing bases, one describes $f$ as an $n \times m$ matrix with coefficients in $A$. For $r \leq \min(m,n)$ consider the ideal in $A$ which is generated by all the $r \times r$ minors of the matrix. It is well known that this ideal $I_{r}(f)$ does not depend on the choice of bases; the largest $r$ for which $I_{r}(f)$ does not vanish is called rk f, the rank of $f$.   

Let
\begin{displaymath}
{\bf F} \quad = \quad 0 \rightarrow F_{s} \rightarrow \cdots \rightarrow F_{i}
\stackrel{d_{i}} \rightarrow F_{i-1} \rightarrow \cdots \rightarrow F_{0}
\end{displaymath}
be a complex of free modules. The by now classical Buchsbaum-Eisenbud
criterion \cite{BE}, as expressed in \cite[Th. 1]{Br}, tells us precisely when this complex is exact. Recall that the grade of an ideal $I$ is the longest length of a regular sequence contained in $I$; this is equal to $\ext^{-}_{A}(A/I,A)$, the smallest degree for which this Ext is $\neq 0$. The grade of the improper ideal is taken to be $\infty$.
\begin{Theorem} Put $f_{i} = \rk F_{i}, i = 0, ..., s$ and $r_{i} = f_{i} -
f_{i+1} + ... \pm f_{s}$. The complex $\bf F$ is acyclic precisely when
$\gr I_{r_{i}}(d_{i}) \geq i$ for all $i$. In this case $r_{i} = \rk d_{i}$.
\end{Theorem}

Instead of minors of a matrix attached to a map between free modules, let
us look at column ideals, i.e. ideals generated by all elements in one particular column. Easy examples show that the grade of such ideals depends
on choice of bases. This motivates the following definition, where we suppose 
that the complex $\bf F$ is exact and minimal in the sense that all boundary
maps in the complex $k \otimes_{A} \bf {F}$ are null. The latter condition is not essential, but it makes the treatment a bit easier.
\begin{Definition} {\rm Let $\bf F$ now be acyclic and minimal. We say that this complex is `stiff' if, regardless of base choice, $\gr \mathfrak c \geq i$ for
every column ideal $\mathfrak c$ belonging to $d_{i}$, $i = 1, ..., s$. We call
the ring $A$ stiff if every such complex over it is stiff.}
\end{Definition}
\begin{Theorem} Every ring of equal characteristic is stiff.  
\end{Theorem}
This result \cite[Th. 1]{SS 4} has some overlap with earlier work of Evans-Griffith \cite{EG} and Hochster-Huneke \cite{HH 1}, \cite{HH 2} as explained in the paper. Our proof uses Big Cohen-Macaulay modules.
\begin{Remark}{\rm An immediate question springs to mind: what have these two 
theorems to do with one another? The Buchsbaum-Eisenbud criterion is by no means obvious. Yet its proof employs techniques and arguments which are rather common 
in this area of local algebra, and are not concerned with characteristic. Is
stiffness really restricted to equal characteristic? The two statements involving grade appear to be in league, but how? For a more detailed discussion, see \cite[section 6]{SS 4}.}
\end{Remark}

What use is stiffness? Well, for a consistent class of rings (if $A$ belongs to 
this class, so does $A/(x)$ for every non zerodivisor $x \in \mathfrak m$) stiffness is controled by first syzygies. In other words, one only needs
to check that $\gr \mathfrak c \geq 1$ for every column ideal $\mathfrak c$
belonging to $d_{1}$ in  each $\bf F$ over every ring in this
class \cite[Prop. 8]{SS 4}. This makes stiffness more accessible, because it is easily seen to mean that $\Ann z = 0$ for every minimal generator $z$ of a syzygy of finite projective dimension. For Gorenstein rings, this turns out to be equivalent with a condition we have already encountered. A little juggling with FID-envelopes yields a proof of \cite[Prop. 10]{SS 4}:
\begin{Proposition} Let $R$ be a Gorenstein ring. Then $\delta(R/\mathfrak b) = 0$ for every unmixed nonnull ideal of zero divisors $\mathfrak b$ if and only if  $\Ann z = 0$ for every minimal generator $z$ of every syzygy of finite projective dimension. 
\end{Proposition}

Tying this together with Corollary 5.3 we obtain \cite[Th. 2]{SS 4} which brings us back to our main theme. 
\begin{Theorem} The Canonical Element Conjecture holds for all rings of 
a particular characteristic and dimension if and only if all such Gorenstein rings are stiff.
\end{Theorem} 
Only the mixed characteristics in dimensions $\geq 4$ remain undecided. 
In all other cases, Gorenstein rings are stiff. However, in equal characteristic   Theorem 6.3 asserts a stronger result. Several obvious questions remain, inviting further exploration.


\smallskip

\begin{thebibliography}{99}

\bibitem{AB}{M. Auslander and R.-O. Buchweitz, The homological theory of maximal Cohen-Macaulay approximations, M\'{e}moire 38, Soc. Math. France (1989), 5-37}.

\bibitem{Ba}{H. Bass, On the ubiquity of Gorenstein rings, Math. Z. 82 (1963),
8-28.}

\bibitem{Be}{P. Berthelot, Alt\'{e}rations de vari\'{e}tes alg\'{e}briques [d'apr\`{e}s A. J. de Jong], S\'{e}m. Bourbaki, Ast\'{e}risque 241, Soc. Math. France 1997, 273-311.}

\bibitem{Br}{W. Bruns, The Evans-Griffith syzygy theorem and Bass numbers,
Proc. Amer. Math. Soc. 115 (1992), 939-946}.

\bibitem{BH}{W. Bruns and J. Herzog, {\it Cohen-Macaulay rings},
Camb. stud. in adv. math. 39, rev. ed., Camb. Univ. Press 1998.}

\bibitem{BE}{D. A. Buchsbaum and D. Eisenbud, What makes a complex exact?, J.
algebra 25 (1973), 259-268}. 

\bibitem{EG}{E. G. Evans and P. Griffith, Order ideals, in {\it Commutative Algebra}, MSRI Publ. Vol. 15, Springer 1989, 212-225}.

\bibitem{Ha}{M. Hashimoto, {\it Auslander-Buchweitz approximations of equivariant modules}, Lond. Math. Soc. Lect. Notes 282, Camb. Univ. Press 2000}.

\bibitem{He}{R. C. Heitmann, The direct summand conjecture in dimension 3, 
Ann. Math. (2) 156 (2002), 695-712}.

\bibitem{Ho 1}{M. Hochster, Contracted ideals from integral extensions of 
regular rings, Nagoya Math. J. 51 (1973), 25-43.}

\bibitem{Ho 2}{M. Hochster, {\it Topics in the homological theory of modules over commutative rings}, Reg. Conf. Ser. Math. 24, Amer. Math. Soc. 1975.}

\bibitem{Ho 3}{M. Hochster, Canonical elements in local cohomology modules and the Direct Summand Conjecture, J. Algebra 84 (1983), 503-553.}  

\bibitem{HH 1}{M. Hochster and C. Huneke, Tight closure, in {\it Commutative 
Algebra}, MSRI Publ. Vol. 15, Springer 1989, 305-324}.

\bibitem{HH 2}{M. Hochster and C. Huneke, Tight closure, invariant theory, and the Brian\c{c}on-Skoda theorem, J. Amer. Math. Soc. 3 (1990), 31-116}.

\bibitem{PS}{C. Peskine and L. Szpiro, Dimension projective finie et cohomologie locale, Publ Math. Inst. Hautes \'{E}tud. Scient. 42 (1973), 77-119}.

\bibitem{Ro 1}{P. Roberts, The Homological Conjectures, in {\it Free Resolutions in Commutative Algebra and Algebraic Geometry}, Research Notes in 
Math. 2, Jones and Bartlett 1992, 121-132.}  

\bibitem{Ro 2}{P. C. Roberts, Recent developments on Serre's multiplicity conjectures: Gabber's proof of the nonnegativity conjecture, Enseignem. Math. 
(2e s\'{e}rie) 44 (1998), 305-324.}

\bibitem{Ro 3}{P. Roberts, Heitmann's proof of the direct summand conjecture in 
dimension 3, arXiv:math.AC/0212073 2002, 1-14}.

\bibitem{Se}{J.-P. Serre, {\it Alg\`{e}bre locale - Multiplicit\'{e}s}, 
Lect. Notes Math. 11, Springer 1965.}                                                                                       
\bibitem{SS 1}{A.-M. Simon and J. R. Strooker, Monomial conjecture and Auslander's $\delta$-invariant, in {\it Commutative algebra and algebraic geometry, ed. F. van Oystaeyen}, Lect. Notes. pure appl. Math. 206, Marcel
Dekker 1999, 265-273}.

\bibitem{SS 2}{A.-M. Simon and J. R. Strooker, Connections between reduced Bass numbers and Auslander's $\delta$-invariant, J. pure appl. alg. 152 (2000),
303-309}.

\bibitem{SS 3}{A.-M. Simon and J. R. Strooker, Reduced Bass numbers, Auslander's $\delta$-invariant and certain homological conjectures, J. reine angew. Math.
551 (2002), 173-218}.

\bibitem {SS 4}{A.-M. Simon and J. R. Strooker, Stiffness of finite free resolutions and the canonical element conjecture, Bull. sciences. math. 127
(2003), 251-260}.

\bibitem{St}{J. R. Strooker, {\it Homological questions in local algebra}, 
Lond. Math. Soc. Lect. Note Ser. 145, Camb. Univ. Press 1990}.

\bibitem{St2}{J. R. Strooker and J. St\"{u}ckrad, Monomial conjecture and complete intersections, Manuscr. Math. 79 (1993), 153-159}.
 
\bibitem{Va}{W. V. Vasconcelos, {\it Divisor theory in module categories}, 
Math. Studies 14, North-Holland, 1973.}


\end{thebibliography}
\end{document}